\documentclass[a4paper,11pt]{article}
\usepackage[english]{babel}
\usepackage{amsmath,amsthm,amsfonts,amssymb}
\usepackage{graphicx}
\usepackage{cite}
  \begin{document}    

\title{A numerical solution of one class of Volterra integral equations of the first kind in terms of the machine arithmetic features                              \thanks{The work was  supported by the RFBR, project  №~15-01-01425-a.}}

\author{Svetlana V. Solodusha, Igor V. Mokry\\
(Melentiev Energy Systems Institute SB RAS, Russia)}
\date{}
\maketitle

\begin{abstract}
 The research is devoted to a numerical solution of the Volterra equations of the first kind that were obtained using the Laplace integral transforms for solving the equation of heat conduction. The paper consists of an introduction and two sections. The first section deals with the calculation of kernels from the respective integral equations at a fixed length of the significand in the floating point representation of a real number. The PASCAL language was used to develop the software for the calculation of kernels, which implements the function of tracking the valid digits of the significand. The test examples illustrate the typical cases of systematic error accumulation. The second section presents the results obtained from the computational algorithms which are based on the product integration method and the midpoint rule.  The results of test calculations are presented to demonstrate the performance of the difference methods.  

\end{abstract}


\section*{Introduction}
\hspace{0.7 cm}

The paper is devoted to the studies on a special class of  Volterra integral equations of the  first kind
  \begin{equation}
 \int\limits_{0}^{t}K_N(t-s)\phi ( s)ds  ={y}(t),\,\,\label{solod31}
\end{equation}
\begin{equation}
 K_N(t-s)= \sum\limits_{q = 1}^{N} {\left( {
- 1} \right)^{q + 1} q^{2}e^{ - \pi ^2 q^2\left( {t - s} \right)}}\label{solod41}.
\end{equation}
The specific feature of kernel $K_N(t-s)$ in \eqref{solod41} lies in that $K_N=0$ in some neighborhood of  zero.
The qualitative theory and numerical methods of solving the Volterra equations of the first kind are dealt with in many studies (for example, \cite{b1986,b2004,v1986,a1999} and the eqreferences given in them). 

The goal of this paper is to consider the application of numerical methods for solving the equations of form \eqref{solod31}, \eqref{solod41} taking into account the mechanisms of error occurrence in the computer calculations.  The work is a continuation of the research started in \cite{s2014, s2015}. 

The Volterra equation of the convolution type \eqref{solod31}, \eqref{solod41} was first obtained in \cite{ya2013}.
The authors of \cite{ya2013} suggest a method of searching for a solution $u(1,t)=\phi(t)$, $t\geq 0$, of an inverse boundary-value problem 
\begin{equation} 
u_{t} = u_{xx} ,\quad x \in \left( {0,1} \right), \quad t \geqslant
0,\label{solod1}
\end{equation}
\begin{equation} 
 u\left( {x,0} \right) = 0, \; u\left( {0,t} \right) = 0, \; 
 u_{x} \left( {0,t} \right) = g_{0}\left( {t} \right)\;
\label{solod2}
\end{equation}
 by reducing \eqref{solod1}, \eqref{solod2} to the Volterra integral equation of convolution type: 
 \begin{equation}
  \int\limits_{0}^{t}K(t-s)\phi ( s)ds  ={y}(t),\,\, 0 \leq s \leq t\leq T,\label{solod3}
\end{equation}
\begin{equation}
 K(t-s)= \sum\limits_{q = 1}^{\infty} {\left( {
- 1} \right)^{q + 1} q^{2}e^{ - \pi ^2 q^2\left( {t - s} \right)}},\; {y}(t)=\frac{1}{2\pi^2}g_{0}\left( {t} \right).\label{solod4}
\end{equation}
Instead of $g_{0}(t) $ we normally know $g_{\delta }(t) $:  
$
{\left\| {g_{\delta}(t)  -
g_{0}}(t) \right\|}_C \leqslant \delta$, $\delta > 0.$

 Problem \eqref{solod1}, \eqref{solod2} plays an important part in the applied problems, including those related to the research into non-stationary thermal processes. Solving the inverse problems is as a rule complicated by the instability of these problems with respect to the initial data errors.
 The application of the methods which employ the Fourier and Laplace transforms in combination with the theory of ill-posed problems found wide application in the construction of stable solutions to the inverse problems of heat conductance. 
 In particular, to solve the problem similar to  \eqref{solod1}, \eqref{solod2},  the authors of \cite{jonas2000} used a stabilizing functional after applying the Fourier transform. In \cite{prud1999} to regularize and assess the convergence of the obtained solutions the authors used a method of conjugate gradients.
In \cite{kolod2010} and \cite{cial2010} the Laplace transform was considered for solving the Cauchy problem. In \cite{monde2003} the Laplace transform was used in a two-dimensional problem. The existing approaches, as a rule, after taking the Laplace transform, apply regularization methods to the obtained equations and then perform the inverse transform.

   Taking into account the ideas from \cite{be2012,ka2011} the authors of  \cite{ya2013}  approximated   $u_x(0,t)$ by the sum $N$ of the first summands: 
 \begin{equation*}
u_x(0,t)= 2 \pi ^2 \sum\limits_{q = 1}^{N} {\left( {
- 1} \right)^{q + 1}q^{2} \int\limits_{ 0}^{t} e^{ - \pi ^2 q^2\left( {t- s} \right)}}\phi(s)ds,
\end{equation*}
where $N$ is positive integer. Then \eqref{solod3},  \eqref{solod4} are reduced to the form  \eqref{solod31},  \eqref{solod41}.
The performance of this approach was discussed in \cite{ya2013}. According to the conclusion made by the author, such a way of solving the inverse problem makes it possible to reduce the initial problem to the Volterra integral equation of the first kind and exclude the components of the operator calculus from the regularization process.

   It is known that the Volterra integral equations of the first kind belong to the class of conditionally-correct problems, and the discretization procedure has a regularizing feature with a regularization parameter, a step of mesh, which is in a certain manner connected to the level of disturbances of the initial data 
$\delta$.

This paper presents an algorithm to numerically solve \eqref{solod31}, \eqref{solod41} at an exactly specified right-hand part.
Note that when solving \eqref{solod31}, \eqref{solod41} we face three types of errors related first of all to the approximation of the initial problem \eqref{solod3}, \eqref{solod4}, secondly to the accuracy of a numerical method, and finally to the computation errors in the machine arithmetic operations with real floating-point numbers. For the research, of greatest interest is the first of the indicated cases. However, to pass to the problem of assessing the parameter $N$ in  \eqref{solod41} it is necessary to develop an algorithm for computation of   $K_{N_{\max}}$, which takes into account the specific features of machine arithmetic  and provides the desired (specified) number of valid digits in the significand.

\section{The specific features of the numerical \\ calculations}
\hspace{0.7 cm}

The computational experiment in \cite{s2014} shows that with an increase in the number of summands in \eqref{solod41} the roots $\lambda^*$ of the equation $K_N(\lambda)=0$ decrease (at the same time monotonisity is observed only separately in even and odd  $N$). The values  $\lambda^*$ will be used further to limit the magnitude of the mesh step $h$ from above, for the value of the mesh function  $K^h_N$ at the first node to be non-zero. As is known, the condition $K_N(0)\neq 0$ is necessary for \eqref{solod31}, \eqref{solod41} to be correct on the pair $(C,C^{1}_{[0,T]})$, where $y(0)=0$, $y'(t)\in C_{[0,T]}$.  

Note, that there can be computational errors in the calculation of $\lambda^*$, they can be related to the application of a fixed mesh in the machine number representation.   

Let us consider the known cases of systematic error accumulation \cite{k1978} which appear in the calculation of kernel values \eqref{solod41}. Use the system of computer algebra Maple10.  Following \cite{m2009} we will include parameter $f$ in the generally accepted representation of the real number. The parameter is equal to the number of valid digits in the significand (starting from the left). Assume that the real number $x=s\cdot M \cdot 10^{-L+p}$   is specified by the set  $(s,M,p,f)$, where  $s\in \{-1,0,+1\}$ is the sign of the number, $M\in\{10^{L-1}, 10^{L-1}+1,...,10^L-1 \}\cup \{0\}$ is significand of the number, 
 $L$ is number of significand positions, $p$ is exponent part of the number.
 
 Let us illustrate the details of the calculations when summing up the numbers with different exponent parts in \eqref{solod41}, using the example.

Example 1.
 Let  $N=50$, $\lambda_0=10^{-3}$, $L\geq 8$.  Take   
\begin{equation*}
x_1=   \sum\limits_{q = 11}^{34} {\left( {
- 1} \right)^{q + 1} q^{2}e^{ - \pi ^2 q^2 {\lambda_0} }},\;
x_2= \sum\limits_{q = 35}^{50} {\left( {
- 1} \right)^{q + 1} q^{2}e^{ - \pi ^2 q^2 {\lambda_0}}}
\end{equation*}
and find  $x_\Sigma=x_1+x_2$.

Assuming 
$$10^{p_\Sigma-f_\Sigma}=10^{p_1-f_1}+10^{p_2-f_2},\; p_{\Sigma}\geq p_1\geq p_2,$$
 according to \cite{m2009},   it is easy to obtain that 
\begin{equation}
f_\Sigma\geq f_s=[f_1- \lg (1 + 10^{-p_1+f_1+p_2-f_2})],\label{ssol5}
\end{equation}
where the symbol $[...]$ means the greatest integer. The last column of the table shows the values of minorant $f_s$, that are calculated using  \eqref{ssol5}. 
 Tab. 1 presents the parameters  $(1,M_1,2,f_1)$, $(1,M_2,-2,f_2)$ and $(1,M_{\Sigma},2,f_{\Sigma})$, which define $x_1$, $x_2$ and $x_\Sigma$ respectively.
\begin{figure}[h!]
\begin{flushright}
\bf{Table 1}
\end{flushright}\centerline{Values $M$  and  $f$  for $x_1$, $x_2$   and $x_\Sigma$.}
\begin{center}
\begin{tabular}{|c|l|c|l|c|l|c|c|}\hline
$L$&    \;\;\;\; \;\; \;\;\;$M_1$&$f_1$ &  \;\;\;\; \;\;\;\;\;  $M_2$&$f_2$&     \;\;\;\;\;\;\;\;\;  $M_\Sigma$&$f_\Sigma$&$f_s$\\  \hline 
8  &18652239 &6 &44981421&6 &18656737 &6&5  \\
9  &186522441 &8 &449814458&7 &186567422 &8&7  \\
10  &1865224455 &9 &4498144699&8 &1865674269 &8 & 8\\
11  &18652244592&11 &44981446726&8 &18656742737 &11 & 10\\
12  &186522445926 &11 &449814466957&10 &186567427373 &11 & 10\\
13  &1865224459248 &12 &4498144669376&10 &1865674273715 &12  & 11\\
14  &18652244592468 &13 &44981446694089&11 &18656742737137 &13  & 12
\end{tabular}
\end{center}
\end{figure}

 The next example illustrates the situation arising in the calculation of the difference between the numbers which have coinciding exponent parts and several high-order digits of the significand.

Example 2. Let $N=50$, $\lambda_0=10^{-3}$, $L\geq 8$. Introduce     
$$
x_3=   \sum\limits_{q = 1}^{10} {\left( {
- 1} \right)^{q + 1} q^{2}e^{ - \pi ^2 q^2 {\lambda_0} }}, \;\;x_4=   \sum\limits_{q = 11}^{50} {\left( {
- 1} \right)^{q + 1} q^{2}e^{ - \pi ^2 q^2 {\lambda_0} }} .
$$
Define $x_\Delta=|x_4|-|x_3|$. 
Suppose that 
$$|M_4-M_3|<10^{L-1}-1, \;p_\Delta\leq p_3=p_4$$
 and, following \cite{{m2009}}  use  the empiric estimate: 
\begin{equation}\label{solod8}
 f_\Delta\geq f_r= \left\{ {\begin{array}{l}
 { \; 0,\;  \text{if}\; \, f_{\min}-L+\lg(\lambda)\leq 0,} \\
 {\left[ f_{\min}-L+\lg(\lambda)\right] ,\;   \text{if}\;\, f_{\min}-L+\lg(\lambda)>0},
  \end{array}}\  \right.
\end{equation}
 where $\lambda=|M_4-M_3|+1,$ $f_{\min}=\min\{f_3,\,f_4\}$.
 
Below are the parameters  $(-1,M_3,2,f_3)$, $(1,M_4,2,f_4)$ and   $(-1,M_\Delta,2,f_\Delta)$, which specify the values  $x_3$,  $x_4$ and $x_\Delta$.
 The estimation from below of $f_r$, obtained using \eqref{solod8}, is given in the last column of Tab. 2. 
 
 \begin{figure}[h!]
\begin{flushright}
\bf{Table 2}
\end{flushright}\centerline{Values  $M$ and  $f$ for $x_3$, $x_4$   and $x_\Delta$.}
\begin{center}
\begin{tabular}{|c|l|c|l|c|l|c|c|}\hline
$L$&    \;\;\;\; \;\; \;\;\;$M_3$&$f_3$ &  \;\;\;\; \;\;\;\;\;  $M_4$&$f_4$&     \;\;\;\;\;\;\;\;\;  $M_\Delta$&$f_\Delta$&$f_r$\\  \hline 
8  &18656743 &7&18656737&6 & 00000006&0  &0\\
9  &186567428 &8&186567424&8 & 000000004&0  &0\\
10  &1865674274 &9 &1865674268&8 & 0000000006&0  &0\\
11  &18656742750 &11 &18656742736&10 & 00000000014&1&0\\
12  &186567427505 &12&186567427372&11 & 000000000133&3 & 1\\
13   &1865674275054 &12&1865674273715&12 &0000000001339&4 & 2\\
14   &18656742750534 &14&18656742737138&13 &00000000013396&4 & 3
\end{tabular}
\end{center}
 \end{figure}

It is obvious that when several high-order digits are set to zero, there appears the number with lower quantity of significant digits in the significand. Fig. 1 illustrates an instantaneous loss of high-order digits, which takes place in this situation. The calculations are made using the software developed by the authors in PASCAL. 
   \begin{figure}[h!]
   \centering
  \includegraphics[scale=0.75]{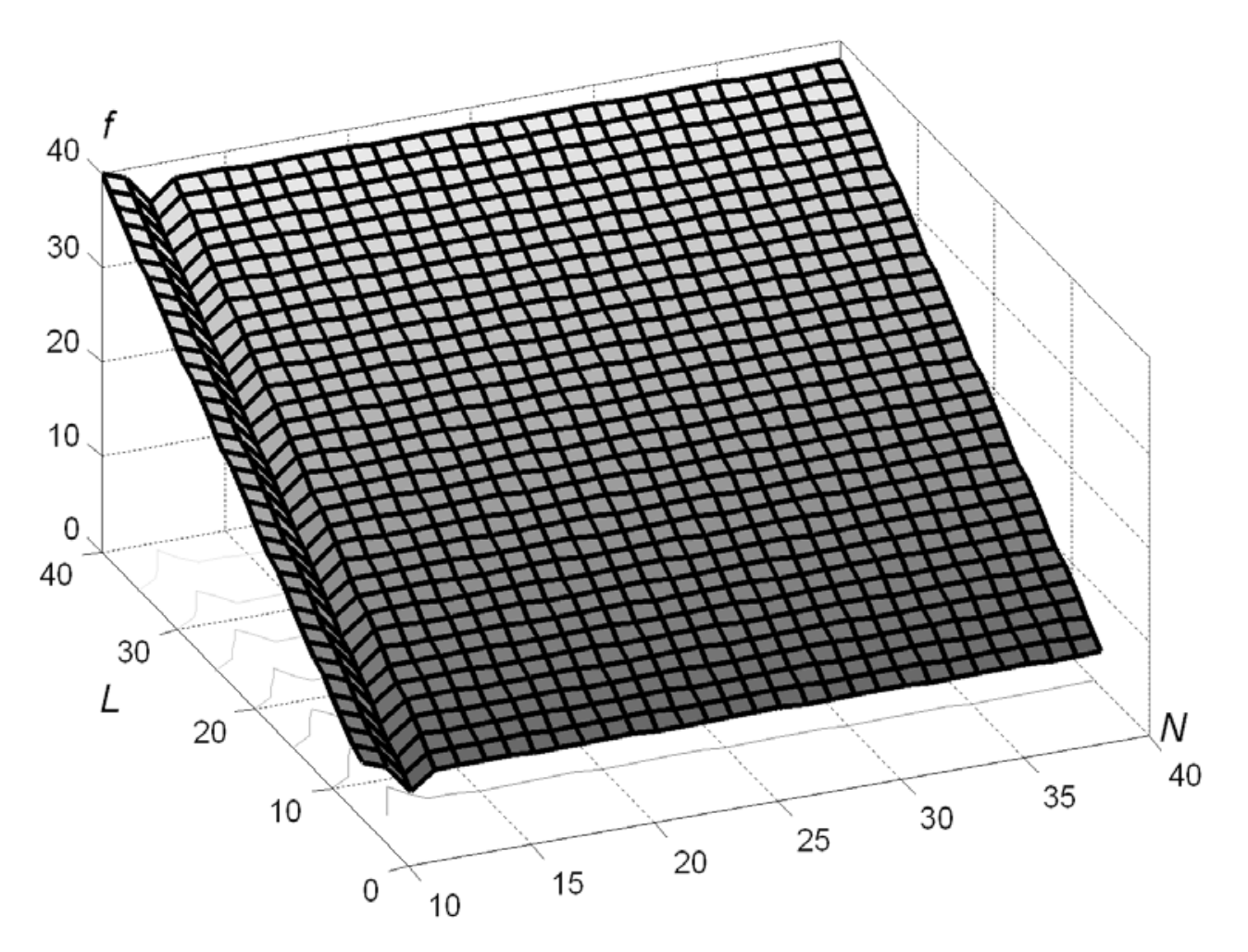}
       \caption{ Instantaneous loss of high-order digits for $K_N(0.001),$ where $ N=12$. }
\end{figure}

 The next section is devoted to the problem of approximately solving \eqref{solod31}, \eqref{solod41} in terms of the specific features of machine arithmetic.   

\section{Results of solving \eqref{solod31}, \eqref{solod41}}
\hspace{0.75 cm}

Let us introduce the uniform meshes of nodes  $t_i=ih,$ $t_{i-\frac{1}{2}}=(i-\frac{1}{2})h,$ $i=\overline{1,n},$  $nh=T$ in $[0,T]$ and, by approximating the integral in \eqref{solod31} by the middle rectangle quadrature and product integration method \cite{linz}, write the corresponding mesh analogues to \eqref{solod31}, \eqref{solod41} as 
 \begin{equation}\label{solod9}
h\sum\limits_{q = 1}^{N}\left( {
- 1} \right)^{q + 1}q^{2}\sum\limits_{j=1}^{i} \,{ {e^{ - \pi ^2 q^2{(j-\frac{1}{2})} h }} \,{\phi}^h_{{i-j+\frac{1}{2}}}}={y}^h_i 
 \end{equation}
 and 
\begin{equation}\label{solod10}
\sum\limits_{q = 1}^{N}\left( {
- 1} \right)^{q + 1}q^{2}\sum\limits_{j=1}^{i} \,\phi^h_{j-\frac{1}{2}}\, {\int\limits_{(j-1)h}^{jh} {e^{ - \pi ^2 q^2{(ih- s)}}}ds} =y^h_i.
 \end{equation}
Designate their solutions by $\check{\phi}^h$ and $\hat{\phi}^h$, respectively. Conduct a numerical experiment with the guaranteed number of valid digits in the significand not less than $f_5=8$.

Example 3.
Assume  $\overline{\phi} \left( {t} \right)$  from \cite{g2010} as the reference function: 
\begin{equation*}
\overline{\phi} \left( {t} \right)=\frac{1-e^{-\frac{t}{\alpha}}}{1-e^{-\frac{1}{\alpha}}}-t, 
 \end{equation*}
 where  $\alpha=10^{-1},\,10^{-2}.$  

Set  in \eqref{solod9}, \eqref{solod10} $N=\overline{2,5};\, 10;\, 15$. Tab. 3, 4 gives the results of numerical calculations in the integration interval $[0,1]$. Here, the notation is as follows: $$
\gamma_1=\log_2 \frac{\Vert \check \varepsilon^{h_1}\Vert_{C_h}}{\Vert \check \varepsilon^{h_2}\Vert_{C_h}}, \; \gamma_2=\log_2 \frac{\Vert \hat \varepsilon^{h_1}\Vert_{C_h}}{\Vert \hat \varepsilon^{h_2}\Vert_{C_h}},
$$ 
where $h_1=2h_2$, $\Vert \check \varepsilon^{h}\Vert_{C_h}$ is maximum absolute difference between the accurate solution and the approximate solution at the nodal points (the approximate solution is obtained using the middle rectangle quadrature); $\Vert \hat \varepsilon^{h}\Vert_{C_h}$ is  maximum absolute difference between the accurate solution and the approximate solution at the nodal points (the approximate solution is obtained using the product integration method); symbol  $*$  means that the error norm is higher than $\max\limits_{0\leq t\leq 1} |\overline{\phi} \left( {t} \right)|$.

\begin{figure}[th]
\begin{flushright}
\bf{Table 3}
\end{flushright}\centerline{Errors in the mesh solution for function  $\overline{\phi}$, where  $\alpha=10^{-1}$.}
\begin{center}
\begin{tabular}{|c|c|c|c|c|c|c|c|c|}\hline
$h$&$||\check\varepsilon||_{C^h}^{N=2}$&$\gamma_1$&$||\hat \varepsilon||_{C^h}^{N=2}$&$\gamma_2$&$||\check\varepsilon||_{C^h}^{N=3}$&$\gamma_1$&$||\hat\varepsilon||_{C^h}^{N=3}$&$\gamma_2$ \\  \hline 
 ${1/256}$  & 0.005001& 2.009  &0.000499& 1.996  &0.003815& 2.002  & 0.001171&1.994 \\ 
 ${1/512}$   & 0.001242& 2.002  &0.000125& 1.998  &0.000952& 2.000  &  0.000294& 1.998\\ 
 ${1/1024}$   & 0.000310&  2.000 &0.000031&1.999   & 0.000238& 2.000  &  0.000074&1.999 \\
 ${1/2048}$   & 0.000078&  1.981 & 0.000008& 1.999  & 0.000059 &1.989 & 0.000018 & 2.001 \\\hline
$h$&$||\check\varepsilon||_{C^h}^{N=4}$&$\gamma_1$&$||\hat\varepsilon||_{C^h}^{N=4}$&$\gamma_2$&$||\check\varepsilon||_{C^h}^{N=5}$&$\gamma_1$&$||\hat\varepsilon||_{C^h}^{N=5}$&$\gamma_2$ \\  \hline 
 ${1/256}$   &  0.065009& 2.107  &0.002056&1.987   &0.025682&2.010&0.003130&1.973 \\ 
 ${1/512}$   &  0.015090& 2.025  &0.000519&1.996   &0.006377&2.002&0.000797& 1.993\\ 
 ${1/1024}$  & 0.003707& 2.006  &0.000129& 1.999  &0.001591&2.000&0.000200& 1.998\\
 ${1/2048}$   & 0.000923&  1.999 &0.000032&  2.000 & 0.000398&1.996&0.000050& 2.000\\\hline
 $h$&$||\check\varepsilon||_{C^h}^{N=10}$&$\gamma_1$&$||\hat \varepsilon||_{C^h}^{N=10}$&$\gamma_2$&$||\check\varepsilon||_{C^h}^{N=15}$&$\gamma_1$&$||\hat\varepsilon||_{C^h}^{N=15}$&$\gamma_2$ \\  \hline 
  ${1/256}$         & $*$& ---  &0.009531& 1.744  &$*$        & ---  & 0.013378&1.394 \\ 
 ${1/512}$   & $*$& ---         &0.002845& 1.922   &0.485248& 2.149  &  0.005092& 1.719\\ 
 ${1/1024}$   & 0.137360&  2.212 &0.000751&1.979   & 0.109395& 2.033  &  0.001547&1.910  \\
 ${1/2048}$   & 0.029650&  2.047 & 0.000190&1.995  & 0.026724 &2.008 &   0.000411 & 1.957  
 \end{tabular}
\end{center}
\end{figure}

\begin{figure}[th]
\begin{flushright}
\bf{Table 4}
\end{flushright}\centerline{Errors in the mesh solution for function  $\overline{\phi}$, where  $\alpha=10^{-2}$.}
\begin{center}
\begin{tabular}{|c|c|c|c|c|c|c|c|c|}\hline
$h$&$||\check\varepsilon||_{C^h}^{N=2}$&$\gamma_1$&$||\hat \varepsilon||_{C^h}^{N=2}$&$\gamma_2$&$||\check\varepsilon||_{C^h}^{N=3}$&$\gamma_1$&$||\hat\varepsilon||_{C^h}^{N=3}$&$\gamma_2$ \\  \hline 
$ 1/256$   & 0.006113&2.009&0.000402&1.995&0.007855&2.004& 0.006159& 1.875\\ 
 $ 1/512$   &0.001518&2.002&0.000101&1.998&0.001958&2.001& 0.001679& 1.939\\ 
$ 1/1024$   & 0.000379&1.999&0.000025&1.992&0.000489&2.000& 0.000438& 1.970\\
$ 1/2048$   & 0.000095&1.985&0.000006& 1.927&0.000122&1.998&0.000112& 1.985\\\hline
$h$&$||\check\varepsilon||_{C^h}^{N=4}$&$\gamma_1$&$||\hat\varepsilon||_{C^h}^{N=4}$&$\gamma_2$&$||\check\varepsilon||_{C^h}^{N=5}$&$\gamma_1$&$||\hat\varepsilon||_{C^h}^{N=5}$&$\gamma_2$ \\  \hline 
$ 1/256$   & 0.080125&2.117&0.014295&1.859 &0.051629&2.022&0.024140& 1.840\\ 
 $ 1/512$   & 0.018474&2.027&0.003940&1.934&0.012716&2.006&0.006744& 1.928\\ 
$ 1/1024$   & 0.004532&2.006&0.001031&1.968&0.003167&2.001&0.001772& 1.966\\
$ 1/2048$   & 0.001128&2.000&0.000264& 1.984&0.000791&2.000&0.000453&1.985\\\hline
 $h$&$||\check\varepsilon||_{C^h}^{N=10}$&$\gamma_1$&$||\hat \varepsilon||_{C^h}^{N=10}$&$\gamma_2$&$||\check\varepsilon||_{C^h}^{N=15}$&$\gamma_1$&$||\hat\varepsilon||_{C^h}^{N=15}$&$\gamma_2$ \\  \hline 
  ${1/256}$  & $*$         & ---  &0.081670& 1.584  &$*$& ---  &  0.114747 &  1.214\\ 
 ${1/512}$   &  $*$        & ---  & 0.027232& 1.849  & $*$ & ---  & 0.049456  &1.637 \\ 
 ${1/1024}$   &0.170835 & 2.232  &0.007556 & 1.945  & 0.2227532   &2.073   & 0.015899 &1.874 \\
 ${1/2048}$   &0.036370 &  2.052 &0.001962 &1.978   & 0.0529318  &2.018 &0.004338 & 1.959 
\end{tabular}
\end{center}
\end{figure}
 
 Tables shows that both finite difference methods have the second order of convergence. 
 Figures 2, and 3 illustrate the behavior of functions  $|\check\varepsilon^h_i|=|\overline{\phi}_{i-\frac{1}{2}}-\check\phi^h_{i-\frac{1}{2}}|$, $|\hat \varepsilon^h_i|=|\overline{\phi}_{i-\frac{1}{2}}-\hat\phi^h_{i-\frac{1}{2}}|$ on a unified mesh with the step   $h=\frac{1}{27}$ at fixed values  $N=2$, $\alpha =10^{-1}$, $T=1$. 
Plot \flqq Line 1\frqq\, corresponds to function $|\hat \varepsilon^h_i|$, plot \flqq Line 2\frqq\, corresponds to function $|\check\varepsilon^h_i|$    and  plot  \flqq Line 3\frqq\, corresponds to  function $| \varepsilon^h_{i_{\min}}|=\min\{|\check\varepsilon^h_i|\, ,|\hat \varepsilon^h_i|\}$.
  \begin{figure}[th]
   \centering
  \includegraphics[scale=0.75]{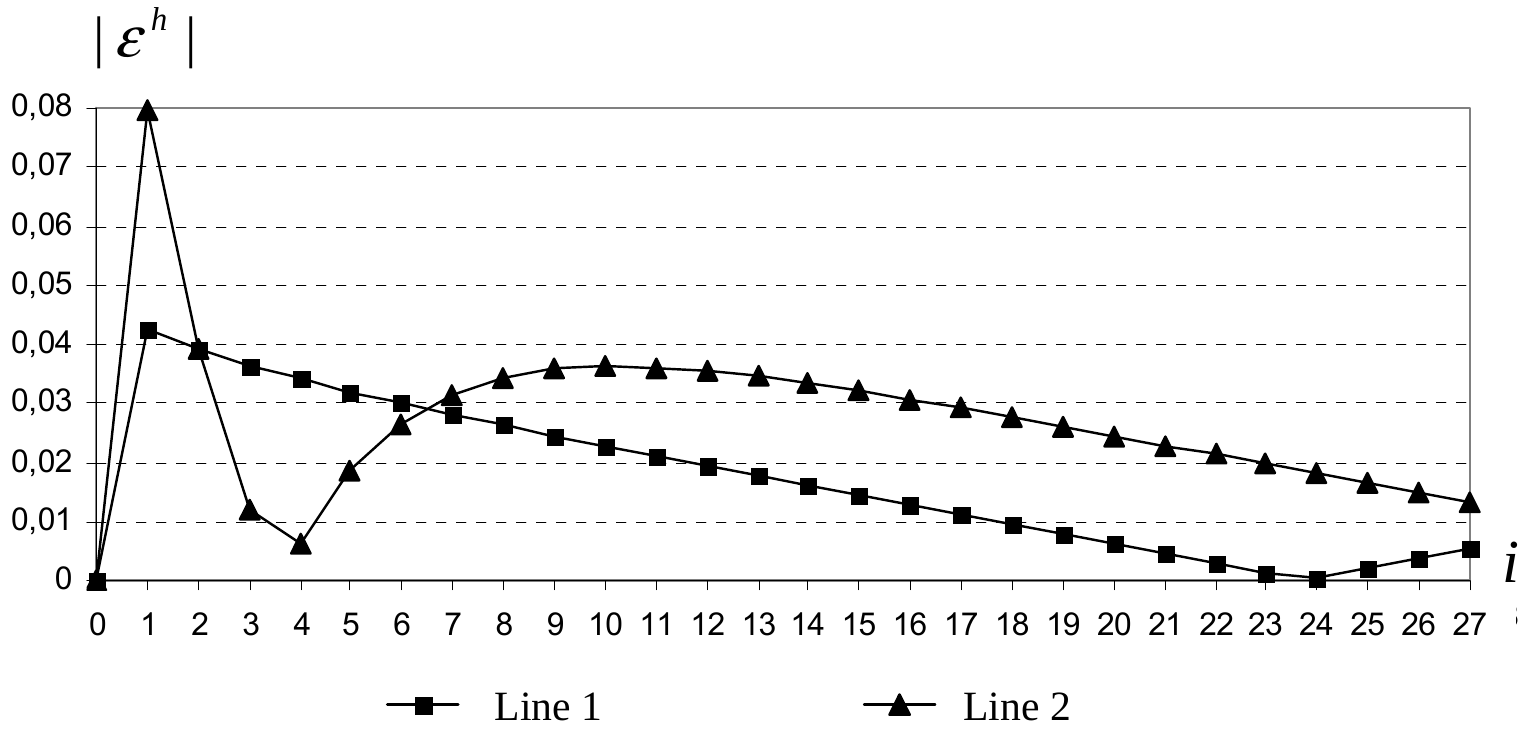} 
     \caption{ Absolute values of errors in mesh solutions at exactly specified initial data.}
\end{figure} 

Figure 2 was obtained at functions $y_i^h$ accurately specified in \eqref{solod9}, \eqref{solod10}. Here, the maximum values of errors made up   $$||\check\varepsilon||_{C^h}^{N=2}=0.0795,\; ||\hat\varepsilon||_{C^h}^{N=2}=0.0424.$$
\begin{figure}[th]
   \centering
  \includegraphics[scale=0.75]{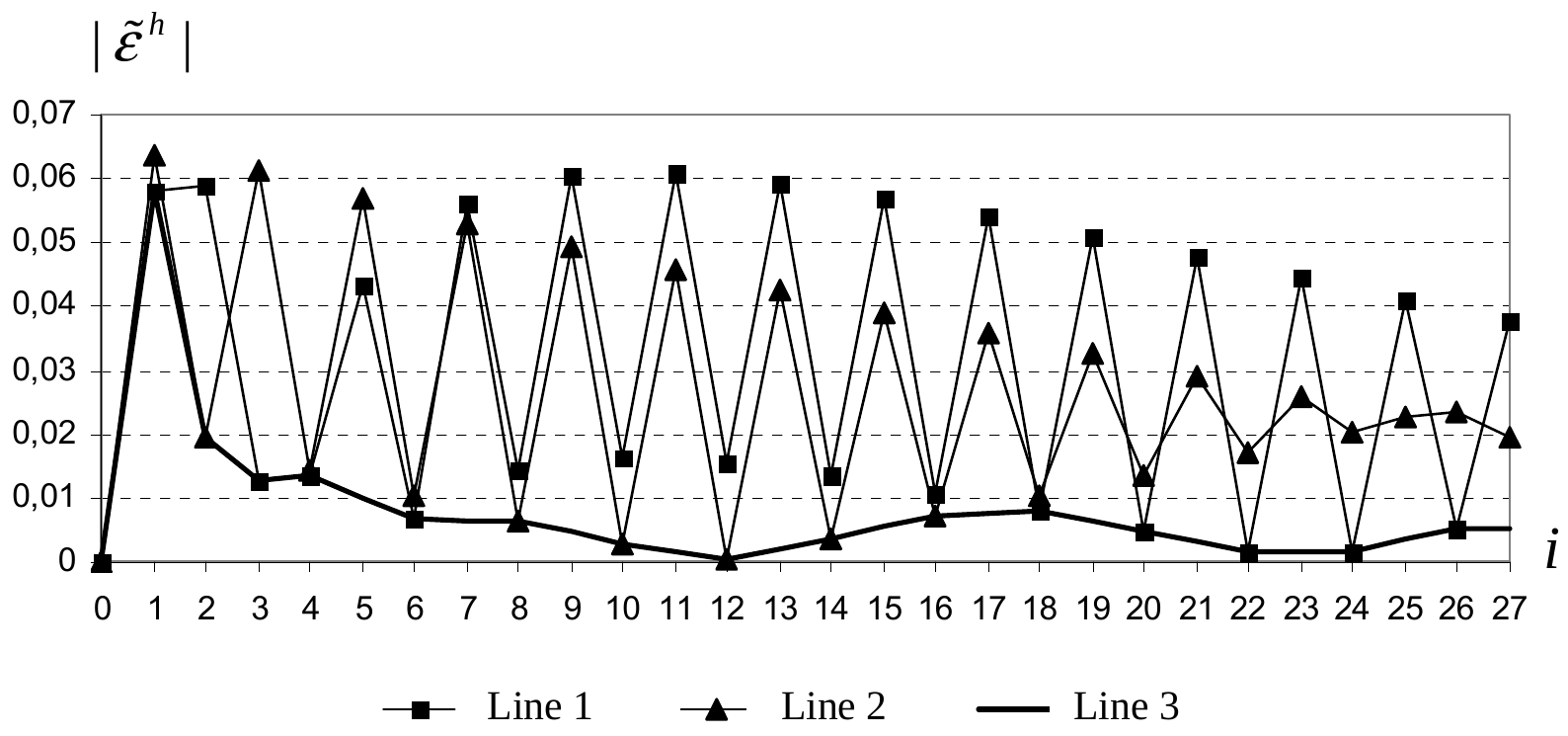}
     \caption{Absolute values of errors in mesh solutions at perturbed initial data. }
\end{figure} 

Figure 3 is obtained at a saw-tooth perturbance of mesh function  $y_i^h$: $$\tilde {y}( t_{i}) = y ( t_{i}) + (  - 1)^{i}\,10^{-3} , \quad i = \overline {1,27,} \,\,\,\,nh = 1.$$ In this case, the maximum values of errors made up
$$||\check\varepsilon||_{C^h}^{N=2}=0.0638,\; ||\hat\varepsilon||_{C^h}^{N=2}=0.0609,\; ||\varepsilon_{\min}||_{C^h}^{N=2}=0.0582.$$ 

  {\bf Remark.  }
   Step  $h$ for the fixed level of the initial data perturbances was chosen using the Fibonacci method with ten trials.
 {\par\addvspace{3mm}}
 The comparison of $\Vert \check \varepsilon^{h}\Vert_{C_h}$ and $\Vert \hat \varepsilon^{h}\Vert_{C_h}$  makes it obvious that the use of the product integration method is more preferable. 
 The computational experiment conducted for the given example at $\alpha=10^{-3}$  shows the convergence for  $h<10^{-3}$. In this connection, further studies are supposed to use the numerical methods of higher order, in particular the third- and fourth-order Runge-Kutta methods. Further, it is planned to construct stability regions of the considered algorithms by the analogy with \cite{bulatov}.

\section*{Conclusion}
\hspace{0.7 cm}
The paper presents the research into the approximate solution of the Volterra integral equation of the first kind of convolution type, which occurs in the inverse boundary-value heat conduction problem, by the second-order finite difference methods. The calculation results obtained using the system Maple 10 are presented. The computational experiment is conducted in terms of the error occurrence mechanism in computer calculations. Typical cases of systematic error accumulation are demonstrated by the test examples. Software for the calculation of kernels, tracking the valid digits in the significand, is developed in PASCAL.


\end{document}